\documentclass[a4paper,12pt]{article}
\usepackage{amssymb,amsmath,amsfonts}
\usepackage{fullpage}
\usepackage[latin1]{inputenc}
\title{Asymptotic tail properties of the distributions in the class of dispersion models}
\author{Alexandre B. Simas$^{a}$\footnote{E-mail: alesimas@impa.br}, Gauss M. Cordeiro$^{b}$\footnote{E-mail: gausscordeiro@uol.com.br} and Saralees Nadarajah$^{c}$\footnote{E-mail: saralees.nadarajah@manchester.ac.uk}\\\\
\centerline{\small{
$^a$Associa\c{c}\~ao Instituto Nacional de Matem\'atica Pura e Aplicada, IMPA,}}\\
\centerline{\small{
Estrada D. Castorina, 110, Jd. Bot\^anico, 22460-320, Rio de Janeiro-RJ, Brasil}}\\
\centerline{\small{$^b$Departamento de Estat\' \i stica e Inform\' atica,
Universidade Federal Rural de Pernambuco,
}}\\
\centerline{\small{
Rua Dom Manoel de Medeiros s/n, Dois Irm\~ aos,
52171-900 Recife-PE, Brasil}}\\
\centerline{\small{$^c$School of Mathematics, University of Manchester}}\\
\centerline{\small{Manchester M13 9PL, UK}}
}
\date{}
\begin{document}
\maketitle
\begin{abstract}
The class of dispersion models introduced by J\o rgensen (1997b) covers many known 
distributions such as the normal, Student t, gamma, inverse Gaussian, hyperbola, von-Mises, 
among others. We study the small dispersion asymptotic (J\o rgensen, 1987b) behavior of the pro\-ba\-bi\-li\-ty density 
functions of dispersion models which satisfy the uniformly convergent saddlepoint approximation. Our 
results extend those obtained by Finner et al. (2008).
\end{abstract}
\section{Introduction}
The class of dispersion models was defined to extend the exponential family models in such a way that
most of their good properties were preserved. As a byproduct of this definition, another class could 
also be defined, the so-called proper dispersion models. The dispersion model, 
defined in the convex support $C$ with position parameter $\mu$ in an 
open interval $\Omega \subset C$ and dispersion parameter 
$\sigma^2>0$, is a family of distributions whose probability density functions (pdfs) 
with respect to some $\sigma$-finite measure 
(usually the Lebesgue or the counting measures) may be written in the form
\begin{equation}\label{dens}
f(y;\mu,\sigma^2) = a(y;\sigma^2)\exp\left\{-\frac{1}{2\sigma^2}d(y;\mu)\right\},\quad y\in C,
\end{equation}
where $a(\cdot;\cdot)\geq 0$ is a suitable function and $d=d(y;\mu)$ is a unit deviance on 
$C\times\Omega$, $\mu\in \Omega$. We denote (\ref{dens}) by the symbol $DM(\mu,\sigma^2)$.
The distribution belongs to the class of proper dispersion models (see J\o rgensen 1997a and 1997b) if its density can be written in form (\ref{dens}) 
with $a(y;\sigma^2)$ being decomposed as $a(y;\sigma^2)= d_1(\sigma^2)d_2(y)$, 
where $d_1(\cdot)$ and $d_2(\cdot)$ are
suitable functions. Further, the exponential dispersion models (see J\o rgensen 1987a and J\o rgensen, 1997b) can be obtained
by taking in (\ref{dens}) $d(y;\mu)=\theta y-b(\theta)+h(y)$, where $\mu = b'(\theta)$ and
$h(y)$ is such that $\theta \mu+h(\mu)=b(\theta)$. The class of proper dispersion models
covers important distributions which are not covered by the exponential dispersion models, 
such as the log-gamma distribution, the Leipnik distribution (Leipnik, 1947 and McCullagh, 
1989) and the reciprocal inverse Gaussian distribution. The von Mises distribution, which also 
belongs to the class of proper dispersion models and does not belong to the exponential dispersion 
models, is particularly useful for the analysis of circular data; see Mardia (1972) and Fisher (1993). 
See J\o rgensen (1997b) for other important examples of proper dispersion models. 
Moreover, the dispersion models have two important general properties. First, the distribution of 
the statistic $D = d(Y;\mu)$ does not depend on $\mu$ when $\sigma^2$ is known, that is, $D$ is a 
pivotal statistic for $\mu$. Second, (\ref{dens}) is an exponential family model with canonical statistic 
$D$ when $\sigma^2$ is known.

\section{Basic Theorems}

First, we introduce some definitions. Let $\Omega \subset C\subset \mathbb{R}$ 
be intervals with $\Omega$ open. A function $d:C\times \Omega\to\mathbb{R}$ is called a 
unit deviance if it satisfies the conditions 
\begin{equation}\label{prop1}
d(y;y) = 0, \quad\forall y\in\Omega\,\,\,\text{and}\,\,\,d(y;\mu) > 0, \quad\forall y\neq\mu.
\end{equation}
Unit deviances may often be constructed on the basis of log-likelihoods. 
In the following, we use the notation $\partial_\mu d$ and $\partial^2_{\mu y} d$ 
to denote the partial derivatives of $d=d(y;\mu)$ with respect to $\mu$ and 
$\mu$ and $y$, respectively, and analogously for higher order derivatives. A unit deviance $d$ is called regular 
if $d(y;\mu)$ is twice continuously differentiable with respect 
to $(y,\mu)$ on $\Omega\times\Omega$ and satisfies
$$\partial^2_{\mu\mu}d(\mu;\mu) > 0, \quad \forall \mu\in\Omega.$$
The unit variance function $V:\Omega\to\mathbb{R}_{+}$ of a regular 
unit deviance is defined by $$V(\mu) = \frac{2}{\partial^2_{\mu\mu}d(\mu;\mu)}.$$

Moreover, it is easy to show that a regular unit deviance $d$ satisfies
\begin{equation}\label{prop3}
\partial_{\mu}d(\mu;\mu) = 0\,\, \hbox{~and~}\,\, \partial_y d(\mu;\mu) = 0,
\end{equation}
for all $\mu\in\Omega$. This holds since (\ref{prop1}) implies that the 
function $d(y;\cdot)$ has a minimum at $y$ and that the function $d(\cdot;\mu)$ 
has a minimum at $\mu$. It can be shown by the chain rule
$$\partial^2_{yy}d(\mu;\mu)=\partial^2_{\mu\mu}d(\mu;\mu)=-\partial^2_{\mu y} d(\mu;\mu),\quad \forall \mu \in\Omega,$$
which gives three equivalent expressions to compute the variance function.

The saddlepoint approximation for a dispersion model with regular unit deviance $d$ is defined 
for $y\in\Omega$ by
\begin{equation}\label{prop4}
f(y;\mu,\sigma^2) \sim \{2\pi\sigma^2 V(y)\}^{-1/2} \exp\left\{-\frac{1}{2\sigma^2} d(y;\mu)\right\}\hbox{~~as~~}\sigma^2\to 0,
\end{equation}
with the notation $a(x)\sim b(x)$ meaning that the ratio $a(x)/b(x)$ goes to $1$ as 
$x\to 0$. The saddlepoint approximation (\ref{prop4}) holds when $\sigma^2\to 0$ for an extensive range of models. 
It is even exact for a few special models, such as the normal and simplex distributions (see Barndorff-Nielsen 
and J\o rgensen, 1991). Moreover, J\o rgensen (1997b) proved that the saddlepoint approximation is uniform on compacts for both proper dispersion models and exponential dispersion models. The asymptotic behavior as $\sigma^2\to 0$ was described 
by J\o rgensen (1987b) as small dispersion asymptotics. The approximation (\ref{prop4}) is equivalent to
$\sigma a(y;\sigma^2) \to \{2 \pi V(y)\}^{-1/2}$ as $\sigma^2\to 0$.
Consider the following theorem whose proof can be found in J\o rgensen (1997b, p. 30).\\
{\bf Theorem 1.} Let $Y \sim DM(\mu_0 + \sigma \mu, \sigma^2)$ be a dispersion model 
with uniformly convergent saddlepoint approximation. Then,
$$\frac{Y-\mu_0}{\sigma}\stackrel{d}{\to} N\{\mu,V(\mu_0)\} \hbox{~~as~~}\sigma^2\to 0,$$
where $\stackrel{d}{\to}$ denotes convergence in distribution. This theorem holds 
for both dispersion and exponential dispersion models. Theorem 1 generalizes 
the well-known result that the density of the Student $t$ random variable with $n$ degrees 
of freedom, say $t_n$, converges to the standard normal density, namely 
$$t_n \stackrel{d}{\to} N(0,1) \hbox{~~as~~}n\to\infty.$$
We can verify this fact by noting that $t_n/\sqrt{n}+\mu$ is a $DM(\mu,1/(n+1))$ 
with unit deviance $d(y;\mu)=\log\{1+(y-\mu)^2\}$. Thus, from the definition of the
variance function we obtain $V(\mu)=1$ and using Theorem 1 it follows
$$\frac{\sqrt{n+1}}{\sqrt{n}} t_n \to N(0,1).$$

Further, in the proof of Theorem 1 given by J\o rgensen (1997b, p. 30), it is shown 
that the pdf of $(Y-\mu_0)/\sigma$, where $Y$ is a 
$DM(\mu_0+\sigma\mu,\sigma^2)$, converges to the normal pdf when $\sigma^2\to 0$, i.e. 
$$f(x;\mu,\sigma^2) \sim \phi(x;\mu,V(\mu_0)),$$
for all $x\in\mathbb{R}$ and $\phi(x;\mu,V(\mu_0))$ stands for the 
pdf of a normal distribution with mean $\mu$ and variance $V(\mu_0)$.
We now give our first results.\\

{\bf Theorem 2.} {\it Let $Y\sim DM(\mu_0+\sigma\mu,\sigma^2)$ be a dispersion 
model with uniformly convergent saddlepoint approximation (in particular, 
the unit deviance is regular), $Z = (Y-\mu_0)/\sigma$, and assume that the unit 
deviance $d=d(y;\mu)$ is of class $C^3$, that is, is continuously three times
differentiable with $\partial^3_{yyy}d(\mu_0;\mu_0)\neq 0$. Thus, if $x_\sigma \in (0,\infty)$ 
for $\sigma \in (0,\infty)$ is such that
$\lim_{\sigma^2\to 0} (x_\sigma-\mu)^3\sigma  = \beta\in [0,\infty]$, then the pdf
$f(x;\mu,\sigma^2)$ of $Z$ satisfies}
$$\lim_{\sigma^2\to 0} \frac{f(x_\sigma;\mu,\sigma^2)}{\phi(x_\sigma;\mu,V(\mu_0))} = \exp\left\{-\frac{\beta}{12}
\partial^3_{yyy} d(\mu_0;\mu_0)\right\}.$$

{\bf Theorem 3.} {\it Let $Y\sim DM(\mu_0+\sigma\mu;\sigma^2)$ be a dispersion model with uniformly convergent 
saddlepoint approximation, $Z = (Y-\mu_0)/\sigma$, and assume that the unit deviance $d$ is of class $C^4$, 
that is, is continuously four times differentiable, with $\partial^4_{yyyy}d(\mu_0;\mu_0)\neq 0$. Suppose further that $\partial^3_{yyy} d(\mu_0+\sigma\mu;\mu_0+\sigma\mu) = 0$ for all $\sigma$ sufficiently small. Thus, 
if $x_\sigma \in (0,\infty)$ for $\sigma \in (0,\infty)$ is such that
$\lim_{\sigma^2\to 0}(x_\sigma-\mu)^4\sigma^2 = \beta\in [0,\infty]$, then the pdf
$f(x;\mu,\sigma^2)$ of $Z$ satisfies}
$$\lim_{\sigma^2\to 0} \frac{f(x_\sigma;\mu,\sigma^2)}{\phi(x_\sigma;\mu,V(\mu_0))} = \exp\left\{-\frac{\beta}{48} 
\partial^4_{yyyy} d(\mu_0;\mu_0)\right\}.$$

More generally,\\

{\bf Theorem 4.} {\it Let $Y\sim DM(\mu_0+\sigma\mu;\sigma^2)$ be a dispersion model with uniformly convergent saddlepoint 
approximation, $Z = (Y-\mu_0)/\sigma$ and assume that the unit deviance $d$ is of class $C^k$, that is, is 
continuously $k$ times differentiable, $k$ a positive integer, with $\partial^k_{y\cdots y}d(\mu_0;\mu_0)\neq 0$, with 
the notation representing $k$ derivatives with respect to $y$. Suppose 
further that $\partial^i_{y\cdots y} d(\mu_0+\sigma\mu;\mu_0+\sigma\mu) = 0$ 
for $\sigma$ sufficiently small and for $2 < i <k$. Thus, 
if $x_\sigma \in (0,\infty)$ for $\sigma \in (0,\infty)$ is such that
$\lim_{\sigma^2\to 0}(x_\sigma-\mu)^k\sigma^{k-2} = \beta\in [0,\infty]$, then the pdf
$f(x;\mu,\sigma^2)$ of $Z$ satisfies}
$$\lim_{\sigma^2\to 0} \frac{f(x_\sigma;\mu,\sigma^2)}{\phi(x_\sigma;\mu,V(\mu_0))} = \exp\left\{\frac{-\beta}{2k!} 
\partial^k_{y\cdots y} d(\mu_0;\mu_0)\right\}.$$

\section{Proof of the theorems}

We shall only prove Theorem 3 since the proofs of Theorems 2 and 4 are entirely analogous.\\
{\bf Proof of Theorem 3:} Let $f_Y(y;\mu,\sigma^2)$ be the pdf of $Y$ and $x_\sigma \in (0,\infty)$ for $\sigma \in (0,\infty)$ be such that
$\lim_{\sigma^2\to 0}(x_\sigma-\mu)^4\sigma^2 = \beta\in [0,\infty]$. Further, the pdf of $Z= (Y-\mu_0)/\sigma$ is
$$f(x;\mu,\sigma^2) = \sigma f_Y(\mu_0+\sigma x;\mu_0+\sigma\mu,\sigma^2),$$
and using (\ref{dens}) we obtain
$$f(x;\mu,\sigma^2) = \sigma a(\mu_0+\sigma x;\sigma^2)\exp\{-(2\sigma^2)^{-1} d(\mu_0+\sigma x;\mu_0+\sigma \mu)\}.$$
Therefore, looking only at $x_\sigma$ from now on, the uniform saddlepoint approximation yields
$$f(x_\sigma;\mu,\sigma^2)\sim \{2\pi V(\mu_0+\sigma x_\sigma)\}^{-1/2} \exp\{-(2\sigma^2)^{-1} d(\mu_0 + \sigma x_\sigma; \mu_0 + \sigma \mu)\},$$
and the continuity of the variance function $V(\cdot)$ gives
\begin{equation}\label{proof1}
f(x_\sigma;\mu,\sigma^2)\sim \{2\pi V(\mu_0)\}^{-1/2}\exp\{-(2\sigma^2)^{-1} d(\mu_0+\sigma x_\sigma;\mu_0+\sigma\mu)\}.
\end{equation}
Thus, from equation (\ref{proof1}) 
$$\frac{f(x_\sigma;\mu,\sigma^2)}{\phi(x_\sigma;\mu,V(\mu_0))} \sim \exp\left\{\frac{(x_\sigma-\mu)^2}{2V(\mu_0)}-(2\sigma^2)^{-1} d(\mu_0+\sigma x_\sigma;\mu_0+\sigma\mu)\right\}.$$
We now focus on the function
$$g_\sigma (x_\sigma)=(2\sigma^2)^{-1} d(\mu_0+\sigma x_\sigma;\mu_0+\sigma\mu).$$
Since the function $d$ is of class $C^4$, the Taylor series expansion up to the 
fourth order in the variable $x_\sigma$ centering on $\mu$ yields 
\begin{eqnarray*}
g_\sigma(x_\sigma) &=& (2\sigma^2)^{-1}\{d(\mu_0+\sigma\mu; \mu_0+\sigma\mu) \\
&+&\sigma (x_\sigma-\mu)\partial_y d(\mu_0 +\sigma\mu;\mu_0+\sigma\mu)+(1/2)\sigma^2(x_\sigma-\mu)^2 \partial^2_{yy} d(\mu_0+\sigma\mu;\mu_0+\sigma\mu)\\
&+& (1/6)\sigma^3(x_\sigma-\mu)^3\partial^3_{yyy} d(\mu_0+\sigma\mu;\mu_0+\sigma\mu)\\
&+& (1/24)\sigma^4(x_\sigma-\mu)^4 \partial^4_{yyyy} d(\mu_0+\sigma\mu;\mu_0+\sigma\mu)+{\cal O}(\sigma^5(x_\sigma-\mu)^5)\}.
\end{eqnarray*}
>From properties (\ref{prop1}) and (\ref{prop3}), we have $d(\mu_0+\sigma\mu;\mu_0+\sigma\mu)=0$ 
and $\partial_y d(\mu_0+\sigma\mu;\mu_0+\sigma\mu)=0$. We also have by assumption 
$\partial^3_{yyy}d(\mu_0+\sigma\mu;\mu_0 +\sigma \mu)=0$, and from the 
definition of the unit variance function $V(\cdot)$, 
it comes $\partial^2_{yy}d(\mu_0+\sigma\mu;\mu_0+\sigma\mu)=2/V(\mu_0+\sigma\mu)$.
Then, we obtain
\begin{eqnarray*}
g_\sigma(x_\sigma) &=& \frac{1}{2}\left\{\frac{(x_\sigma-\mu)^2}{V(\mu_0+\sigma\mu)}\right. \\
&+&\left.\frac{1}{24}\sigma^2(x_\sigma-\mu)^4\partial^4_{yyyy}d(\mu_0+\sigma\mu;\mu_0+\sigma\mu) 
+ {\cal O}(\sigma^3(x_\sigma-\mu)^5)\right\}.
\end{eqnarray*}
Hence, we have
$$\frac{f(x_\sigma;\mu,\sigma^2)}{\phi(x_\sigma;\mu,V(\mu_0))}\sim \exp\left\{-\frac{1}{48}\sigma^2(x_\sigma-\mu)^4 
\partial^4_{yyyy} d(\mu_0+\sigma\mu;\mu_0+\sigma\mu)+{\cal O}(\sigma^3(x_\sigma-\mu)^5)\right\}.$$

Furthermore, the continuity of the fourth derivative together with the fact that 
$\sigma(x_\sigma-\mu)^4 \to \beta$ as $\sigma^2 \to 0$ implies that
\begin{equation}\label{proof2}
\sigma^2 (x_\sigma-\mu)^4 \partial^4_{yyyy}d(\mu_0+\sigma\mu;\mu_0+\sigma\mu)\to \beta \partial^4_{yyyy}d(\mu_0;\mu_0).
\end{equation}
Moreover, the hypothesis that $\sigma^2 (x_\sigma-\mu)^4 \to \beta$ as $\sigma^2\to 0$ also implies that
$${\cal O}(\sigma^3(x_\sigma-\mu)^5) = o(1).$$
Hence, from the last equation and (\ref{proof2}), we obtain
$$\frac{f(x_\sigma;\mu,\sigma^2)}{\phi(x_\sigma;\mu,V(\mu_0))} \sim \exp\left\{-\frac{\beta}{48}\partial^4_{yyyy} d(\mu_0;\mu_0)\right\},$$
which concludes the proof.\\\\

%We now state a corollary whose proof follows the same lines of that from Corollary 1.1 from Finner et al. (2008).\\
%{\bf Corollary:} (Mill's Ratio for $(x_\sigma\to\infty)$. {\it If $\lim_{\sigma^2\to 0} (x_\sigma-\mu)^3\sigma  = \beta\in [0,\infty)$ (or
%$\lim_{\sigma^2\to 0}(x_\sigma-\mu)^4 \sigma^2  = \beta\in [0,\infty)$ depending on the same conditions of Theorem 2 and 3) and
%$\lim_{\sigma^2\to 0} x_\sigma = \infty$, then}
%$$\frac{F_\sigma(-x_\sigma)}{f_\sigma(x_\sigma)} \sim \frac{1}{x_\sigma} \hbox{~as~} \sigma\to 0.$$

\section{Special cases}
In order to apply Theorems 2, 3 and 4 to a special distribution, we write the distribution in the 
form (\ref{dens}) and compute the third derivative, check if it is zero, or not, if so, compute the 
fourth derivative, check if it is zero and so on, and then apply the results. We now provide some 
special cases.\\

\textbf{Student $t$ distribution.}\\

The pdf of the Student $t$ random variable with $n$ degrees of freedom, $t_n$, is given by
$$f(x;n) = \frac{1}{\sqrt{n} B(n/2,1/2)}\left(1+\frac{x^2}{n}\right)^{-\frac{n+1}{2}},$$
where $B(\cdot,\cdot)$ is the beta function. Then, $t_n/\sqrt{n}+\mu_0$ is a 
$DM(\mu_0,1/(n+1))$. A straightforward computation shows that the unit deviance 
function is $d(y;\mu_0)=\log\{1+(y-\mu_0)^2\}$ which implies that 
$\partial_{yy}^2 d(\mu_0;\mu_0) = 2$ and $V(\mu_0)=1$ for all $\mu_0$. We obtain 
$\partial_{yyy}^3 d(\mu_0;\mu_0)=0$ and $\partial_{yyyy}^4 d(\mu_0;\mu_0) = -12$. 
Hence, taking $\mu_0 = 0$ and by Theorem 3, the pdf of 
$(\sqrt{n+1}/\sqrt{n})t_n$, and thus the pdf $f_n(x)$ of $t_n$ 
obeys $f_n(x_\sigma)/\phi(x_\sigma;0,1) \to \exp(\beta/4)$ as $\sigma\to 0$, 
which is the result given in Theorem 1.1 of Finner et al. (2008).\\

\textbf{Generalized Student $t$ distribution}.\\

The pdf of the generalized Student $t$ random variable $X_{r,s}$ with parameters $r$ and 
$s$ is given by
$$f(x;s,r)=\frac{1}{\sqrt{s}B(r/2,1/2)}\left(1+\frac{x^2}{s}\right)^{-\frac{r+1}{2}},$$
with $s,r>0$. It includes Student $t$ distribution with $n$ degrees of freedom when $s=r=n$. For this distribution we 
have two possibilities: first, $s$ is known and fixed, then $X_{r,s}+\mu_0$ is a dispersion model 
with dispersion parameter $\sigma^2 = 1/(r+1)$ and unit deviance $d(y;\mu_0)=\log\{1+(y-\mu_0)^2/s\}$. 
We have $\partial_{yy}^2 d(\mu_0;\mu_0) = 2/s$ which implies that $V(\mu_0) = s$ for all $\mu_0$. 
The derivatives are hence $\partial_{yyy}^3 d(\mu_0;\mu_0)=0$ and $\partial_{yyyy}^4 d(\mu_0;\mu_0)=\frac{-12}{s^2}$
and, by Theorem 3, imply that the pdf $f_r(x)$ of $\sqrt{r+1}(X_r-\mu_0)$ satisfies 
$f_r(x_\sigma)/\phi(x_\sigma;0,s) \to \exp\{\beta/(4s^2)\}$ as $\sigma\to 0$, which is equivalent 
to $r\to\infty$.\\
We now consider the second possibility when $s$ depends on $r$ and, also, $s(r)/r \to a>0$ as 
$r \to\infty$. Hence, $X_{r,s}/\sqrt{s}+\mu_0$ is a $DM(\mu_0,1/(r+1))$ with unit deviance 
$d(y;\mu_0)=\log\{1+(y-\mu_0)^2\}$. Hence, $V(\mu_0) = 1$ for all $\mu_0$, 
$\partial_{yyy}^3 d(\mu_0;\mu_0)=0$ and $\partial_{yyyy}^4 d(\mu_0;\mu_0)=-12$ 
as discussed before for the Student $t_n$ distribution. Taking $\mu_0=0$, the density of 
$\sqrt{(r+1)/s}X_{r,s}$, and then the pdf $g_{r,s}(x)$ of $a^{-1} X_{r,s}$ verifies 
$g_{r,s}(x_\sigma)/\phi(x_\sigma;0,1) \to \exp(\beta/4)$ as $\sigma^2\to 0$, 
or equivalently, as $r \to \infty$. Moreover, the pdf $f_{r,s}(x)$ of $X_{r,s}$ 
satisfies $f_{r,s}(a x_\sigma)/\phi(x_\sigma;0,1)\to a^{-1}\exp(\beta/4)$ as $r\to\infty$.\\

\textbf{Gamma distribution}\\

The pdf of a gamma random variable $X_\lambda$ with mean $\mu_0$ and precision parameter 
$\lambda$ is given by
$$f(y;\mu_0,\lambda) = \frac{\lambda^\lambda e^{-\lambda}}{y \Gamma(\lambda)} \exp\left\{-\lambda\left[ \frac{y}{\mu_0}-\log\left(\frac{y}{\mu_0}\right)-1\right]\right\},$$
which is obviously a $DM(\mu_0,1/\lambda)$ with unit deviance 
$d(y;\mu_0)=2\{y/\mu_0-\log(y/\mu_0)-1\}$. We have 
$\partial_{yy}^2 d(\mu_0;\mu_0) =2/\mu_0^2$ and $V(\mu_0) = \mu_0^2$. 
Further, $\partial_{yyy}^3 d(\mu_0;\mu_0)=-4/\mu_0^3$ and using Theorem 2 it follows 
that the density $f_\lambda(x)$ of $\sqrt{\lambda}(X_\lambda-\mu_0)$ 
satisfies $f_{\lambda}(x_\sigma)/\phi(x_\sigma;0,\mu_0^2) \to \exp\{\beta/(3\mu_0^3)\}$ as 
$\sigma^2\to 0$, or equivalently, as $\lambda\to\infty$.\\

\textbf{Reciprocal gamma distribution}\\

The pdf of the reciprocal gamma distribution with parameters $\mu_0$ and $\lambda>0$ is given by
$$f(y;\mu_0,\lambda)=\frac{\lambda^\lambda e^{-\lambda}}{y\Gamma(\lambda)}\exp\left\{\frac{-\lambda\mu_0}{y} + \lambda\log\left(\frac{\mu_0}{y}\right)+\lambda\right\},$$
This distribution is also known as inverse gamma distribution and is a $DM(\mu_0,1/\lambda)$ with
unit deviance $d(y;\mu_0) = 2(\mu_0/y-\log(\mu_0/y)-1)$. We have $\partial_{yy}^2 d(\mu_0;\mu_0) =2/\mu_0^2$ (such as for the gamma 
distribution), and then the unit variance is $V(\mu_0) = \mu_0^2$. Further, $\partial_{yyy}^3 d(\mu_0;\mu_0)=-8/\mu_0^3$ and 
then, by Theorem 3, the density $f_\lambda(x)$ of $\sqrt{\lambda}(X_\lambda-\mu_0)$ satisfies $f_{\lambda}(x_\sigma)/\phi(x_\sigma;0,\mu_0^2)\to\exp\{2\beta/(3\mu_0^3)\}$ as $\lambda\to \infty$.\\

\textbf{Log-Gamma distribution}\\

The pdf of the log-gamma distribution with parameters $\mu_0$ and $\lambda>0$ is given by
$$f(y;\mu_0,\lambda)=\frac{\lambda^\lambda}{\Gamma(\lambda)}\exp\left\{\lambda(y-\mu_0-\exp\{y-\mu_0\})\right\}.$$
This distribution is a $DM(\mu_0,1/\lambda)$ with unit deviance $d(y;\mu_0) = 2(-y+\mu_0+\exp\{y-\mu_0\}-1)$. 
We have $\partial_{yy}^2 d(\mu_0;\mu_0) =2$ which leads to $V(\mu_0) = 1$. We also have $\partial_{yyy}^3 d(\mu_0;\mu_0)=2$ 
and then, by Theorem 2, the density $f_\lambda(x)$ of $\sqrt{\lambda}(X_\lambda-\mu_0)$ satisfies $f_{\lambda}(x_\sigma)/\phi(x_\sigma;0,\mu_0^2)\to\exp\{-\beta/6\}$ as $\lambda\to \infty$.\\

\textbf{Generalized hyperbolic secant distribution}\\

The pdf of the generalized hyperbolic secant (GHS) distribution is given by
$$f(y;\mu_0,\lambda)=\lambda c^\ast(\lambda y; \lambda)\exp\{\lambda[y\arctan \mu_0 - 0.5\log(1+\mu_0^2)]\},$$
where 
$$c^\ast(z;\lambda) = \frac{2^{\lambda-2}|\Gamma(\lambda/2 + iz/2)|^2}{\pi\Gamma(\lambda)},$$
and $\mu_0$ and $\lambda$ are the mean and the precision parameter, respectively. 
Hence, if $X_\lambda$ follows a GHS distribution, then it is a $DM(\mu_0,1/\lambda)$
for which the unit deviance is 
$$d(y;\mu_0)=2y(\arctan y-\arctan\mu_0)+\log\{(1+\mu_0^2)/(1+y^2)\}.$$ 
Hence, $\partial_y^2 d(\mu_0;\mu_0)=2/(1+\mu_0^2),$ which yields $V(\mu_0)=1+\mu_0^2$.
Further, $\partial_{yyy}^3 d(\mu_0;\mu_0)=-4\mu_0/(1+\mu_0^2)^2$ which implies by 
Theorem 2 that the pdf $f_\lambda(x)$ of $\sqrt{\lambda}(X_\lambda-\mu_0)$ 
satisfies $f_{\lambda}(x_\sigma)/\phi(x_\sigma;0,1+\mu_0^2) \to \exp\{\beta\mu_0/(3(1+\mu_0^2)^2)\}$ 
as $\sigma^2\to 0$, or equivalently, as $\lambda\to\infty$.\\

\textbf{Inverse Gaussian distribution}\\

The pdf is
$$f(y;\mu_0,\sigma^2) = \frac{1}{\sqrt{2\pi\sigma^2}}y^{3/2} \exp\left\{-\frac{(y-\mu_0)^2}{2\sigma^2 y\mu_0^2}\right\},$$
which is a $DM(\mu_0,\sigma^2)$ with unit deviance $d(y;\mu_0)=(y-\mu_0)^2/(y\mu_0^2)$. We have  
$\partial_{yy}^2 d(\mu_0;\mu_0)=2/\mu_0^3$ and then the variance function is $V(\mu_0)=\mu_0^3$. 
We obtain $\partial_{yyy}^3 d(\mu_0;\mu_0)= -6\mu_0^{-4}$. This justifies, by Theorem 2, that the density 
$f(x)$ of $(X-\mu_0)/\sigma$ when $\sigma^2 \to 0$ satisfies $f(x_\sigma)/\phi(x_\sigma;0,\mu_0^3)\to \exp\{\beta/(2\mu_0^4)\}$.\\

\textbf{Reciprocal inverse Gaussian distribution}\\

The pdf is
$$f(y;\mu_0,\lambda)=\sqrt{\frac{\lambda}{2\pi}}y^{-1/2} \exp\left\{-\frac{\lambda(y-\mu_0)^2}{2y}\right\},$$
which is a $DM(\mu_0,1/\lambda)$ with unit deviance $d(y;\mu_0)=(y-\mu_0)^2/y$, and then $\partial_{yy}^2 d(\mu_0;\mu_0)=2/\mu_0,$
yielding the variance function $V(\mu_0)=\mu_0$. We have $\partial_{yyy}^3 d(\mu_0;\mu_0)= -6\mu_0^{-2}$ and 
the following result, due to Theorem 2, for the density $f(x)$ of $(X-\mu_0)/\sigma$ when $\sigma^2\to 0$ holds: 
$f(x_\sigma)/\phi(x_\sigma;0,\mu_0) \to \exp\{\beta/(2\mu_0^2)\}$.\\

\textbf{Hyperbola distribution}\\

The pdf is 
$$f(y;\mu_0,\lambda)=\frac{e^{-\lambda}}{2 K_0(\lambda)}y^{-1} \exp \left\{-\frac{\lambda(y-\mu_0)^2}{2y\mu_0}\right\},$$
where $K_0$ is a Bessel function and $\mu_0$ and $\lambda$ are parameters. It is a 
$DM(\mu_0,1/\lambda)$, $\sigma^2=1/\lambda$, with unit deviance $d(y;\mu_0)=(y-\mu_0)^2/(y\mu_0)$, 
yielding $\partial_y^2 d(\mu_0;\mu_0)=2/\mu_0^2$, and then va\-riance 
function $V(\mu_0)=\mu_0^2$. Further, $\partial_{yyy}^3 d(\mu_0;\mu_0)= -6\mu_0^{-3}$ which leads, by Theorem 2, to the 
following result for the density $f(x)$ of $\sqrt{\lambda}(X-\mu_0)$: $f(x_\sigma)/\phi(x_\sigma;0,\mu_0^2) \to \exp\{\beta/(2\mu_0^3)\}$ 
when $\sigma^2\to 0$.\\

\textbf{Hyperbolic distribution}\\

The pdf for $y \in \mathbb{R}$ is
$$f(y;\mu_0,\lambda) = \frac{1}{2 a K_1(\lambda)} \exp\left( -\lambda\left[a\{1+(y-\mu_0)^2\}^{1/2}-b(y-\mu_0)\right]\right),$$
where $K_1$ is a Bessel function, $a^2=1+b^2$ and $\mu_0$ and $\lambda$ are parameters. This distribution 
is a $DM(\mu_0,1/\lambda)$ for $b$ fixed, $\sigma^2=1/\lambda$, with unit deviance $d(y;\mu_0)=a \{1+(y-\mu_0)^2\}^{1/2}-b(y-\mu_0)-a$
and then $\partial_{yy}^2 d(\mu_0;\mu_0)=a$ and $V(\mu_0)=2/a$. We have $\partial_{yyy}^3 d(\mu_0;\mu_0)=0$ and 
$\partial_{yyyy}^4 d(\mu_0;\mu_0)=-3a$. As $\sigma\to 0$, Theorem 3 implies that the pdf $f(x)$ of $\sqrt{\lambda}(X-\mu_0)$ obeys 
$f(x_\sigma)/\phi(x_\sigma;0,\mu_0) \to \exp(a\beta/4)$.\\

\textbf{Simplex distribution}\\

The pdf of the simplex distribution with parameters $\mu_0 \in (0,1)$ and $\sigma^2>0$ is 
$$f(y;\mu_0,\sigma^2) = \frac{1}{\sqrt{2\pi\sigma^2 \{y(1-y)\}^3}}\exp\left\{-\frac{(y-\mu_0)^2}{2\sigma^2 y(1-y)\mu_0^2(1-\mu_0)^2}\right\}.$$
This distribution is an example of a proper dispersion model suitable for proportions and was introduced by Barndorff-Nielsen 
and J\o rgensen (1991). This distribution is a $DM(\mu_0,\sigma^2)$ with unit deviance $d(y;\mu_0)=(y-\mu_0)^2/\{y(1-y)\mu_0^2(1-\mu_0)^2\}$. We have $\partial_y^2 d(\mu_0;\mu_0)=2/\{\mu_0^3(1-\mu_0)^3\}$ and $V(\mu_0)=\mu_0^3(1-\mu_0)^3$. We obtain $\partial_{yyy}^3 d(\mu_0;\mu_0)=6(2\mu_0-1)/\{\mu_0(1-\mu_0)\}^4$ and then, by Theorem 2, the density $f(x)$ of $(X-\mu_0)/\sigma$ 
when $\sigma^2\to 0$ satisfies $f(x_\sigma)/\phi(x_\sigma;0,V(\mu_0)) \to \exp[-\beta(\mu_0-0.5)/\{\mu_0(1-\mu_0)\}^4]$.\\

\textbf{von-Mises distribution}\\

The pdf of the von Mises distribution may be expressed as
$$f(y;\mu_0,\sigma^2)=\frac{e^{\sigma^{-2}}}{2\pi I_0(\sigma^{-2})} \exp\left\{-\frac{1-\cos(y-\mu_0)}{\sigma^2}\right\},$$
for $ 0 \le y \le 2 \pi$, $\mu_0 \in [0,2\pi)$ and $\sigma^2>0$, where $I_0(\cdot)$ denotes the modified 
Bessel function. This distribution is a $DM(\mu_0,\sigma^2)$ with unit deviance $d(y;\mu_0)=2\{1-\cos(y-\mu_0)\}$.
We have $\partial_{yy}^2 d(\mu_0;\mu_0)=2$ and $V(\mu_0)=1$. It follows $\partial_{yyy}^3 d(\mu_0;\mu_0)=0$, 
$\partial_{yyyy}^4 d(\mu_0;\mu_0)=-2$ and then, by Theorem 3, the density $f(x)$ of 
$(X-\mu_0)/\sigma$ satisfies $f(x_\sigma)/\phi(x_\sigma;0,1) \to \exp(\beta/24)$ when $\sigma^2\to 0$.\\

\textbf{Leipnik distribution}\\

The pdf of the Leipnik distribution for $y,\mu_0 \in (-1,1)$ and $\lambda > 0$ is
$$f(y;\mu_0,\lambda) =\frac{(1-y^2)^{-1/2}}{B\left(\frac{\lambda+1}{2},\frac{1}{2}\right)}\exp\left\{-\frac{\lambda}{2} 
\log\left(\frac{1-2y\mu_0+\mu_0^2}{1-y^2}\right)\right\},$$
where $B(\cdot,\cdot)$ denotes the beta function. Leipnik (1947) derived this as a `smoothed' 
approximation to the distribution for the circular serial correlation coefficient for a sample of size 
$\lambda$. McCullagh (1989) rediscovered the distribution as a kind of noncentral version of the symmetric 
beta family, and noted a connection with Brownian motion. This distribution is a $DM(\mu_0,1/\lambda)$, 
$\sigma^2=1/\lambda$, with unit deviance $d(y;\mu_0)=\log\left(\frac{1-2y\mu_0+\mu_0^2}{1-y^2}\right)$. 
Noting that $\partial_{yy}^2 d(\mu_0;\mu_0)=2/(1-\mu_0^2)$, the variance function reduces to $V(\mu_0)=1-\mu_0^2$. We have $\partial_{yyy}^3 d(\mu_0;\mu_0)= 12\mu_0/\{(1-\mu_0^2)^2\}$, 
and therefore, by Theorem 2, the density $f(x)$ of $\sqrt{\lambda}(X-\mu_0)$ obeys
$f(x_\sigma)/\phi(x_\sigma,0,1-\mu_0^2) \to \exp\{-\beta\mu_0/(1-\mu_0^2)^2\}$ when $\lambda\to\infty$.\\

\textbf{Transformed Leipnik distribution}\\

The pdf of the transformed Leipnik distribution for $y,\mu_0 \in (0,1)$ and $\lambda > 0$ is
$$f(y;\mu_0,\lambda)=\frac{\{y(1-y)\}^{-1/2}}{B((\lambda+1)/2,1/2)}\exp\left[-\frac{\lambda}{2} \log\left\{1+\frac{(y-\mu_0)^2}{y(1-y)}\right\}\right],$$
which is a $DM(\mu_0,1/\lambda)$, $\sigma^2=1/\lambda$, with unit deviance 
$d(y;\mu_0)=\log\left\{1+\frac{(y-\mu_0)^2}{y(1-y)}\right\}$. Also, $\partial_{yy}^2 d(\mu_0;\mu_0)=2/\{\mu_0(1-\mu_0)\}$ and 
then $V(\mu_0)=\mu_0(1-\mu_0)$. Moreover, $\partial_{yyy}^3 d(\mu_0;\mu_0)=6(2\mu_0-1)/\{\mu_0(1-\mu_0)\}^2$. 
Theorem 2 shows that the density $f(x)$ of $\sqrt{\lambda}(X-\mu_0)$ implies that $f(x_\sigma)/\phi(x_\sigma,0,\mu_0(1-\mu_0))
\to\exp\{-\beta(2\mu_0-1)/[12\{\mu_0(1-\mu_0)\}^2]\}$ as $\lambda\to\infty$.\\

\textbf{Modified three parameter generalized inverse Gaussian}\\

The pdf for the generalized inverse Gaussian with parameters $\mu_0,\lambda>0$ and 
$a \in [-1,1]$ is given by
$$f(y;\mu_0,\lambda)=\frac{\left(\frac{1+a}{1-a}\right)^{\lambda a/2} e^{-\lambda/2}}{2y K_{\lambda a}(\lambda\sqrt{1-a^2})}\exp\left\{-\frac{\lambda}{2} d(y;\mu_0)\right\},$$
where $K_{\lambda a}(\cdot)$ denotes the modified Bessel function of third kind with 
index $\lambda a \in \mathbb{R}$. This distribution was considered by J\o rgensen (1982, 1997b). 
This distribution contains the gamma distribution when $a=1$ and the reciprocal gamma distribution 
when $a=-1$. This distribution is a $DM(\mu_0,1/\lambda)$ for fixed $a$, $\sigma^2=1/\lambda$, 
and the unit deviance $d(y;\mu_0)$ reduces to 
$$d(y;\mu_0) = 2 a \log\left(\frac{\mu_0}{y}\right)+\frac{y}{\mu_0}(1+a)+\frac{\mu_0}{y}(1-a)-2.$$
We have $\partial_{yy}^2 d(\mu_0;\mu_0)=2/(\mu_0^2)$ and $V(\mu) = \mu_0^2$. Moreover, 
$\partial_{yyy}^3 d(\mu_0;\mu_0) = -(6-2a)/\mu_0^3$. Therefore, by Theorem 2, the density $f(x)$ of 
$\sqrt{\lambda}(X-\mu_0)$ satisfies $f(x_\sigma)/\phi(x_\sigma,0,\mu_0^2)\to \exp\{(3-a)\beta/(6\mu_0^3)\}$ as $\lambda\to\infty$.\\

\section{Conclusions}

Under the assumption of uniform saddlepoint approximation, dispersion models are approximately normal for small values of the
dispersion parameter $\sigma^2$, the so-called small dispersion asymptotics. This fact with a few extra regularity conditions, 
lead us to gene\-ralize Finner's (2008) result to a much wider class of distributions. Moreover, in a much more flexible 
fashion, since we allow the derivatives of the unit deviance function to vanish any finite amount of times, 
and the result of Theorem 4 still holds. Finally, several special distributions were considered and 
studied in detail.

\end{document}